\documentclass{article}
\usepackage{amsmath}
\input amssymb.sty
\newtheorem{prop}{Proposition}[section]
\newtheorem{lem}[prop]{Lemma}

\def\endproof{\hfill $\Box\qquad$\endtrivlist} 
\newtheorem{pr}[prop]{Property}
\newtheorem{example}[prop]{Example}
\newtheorem{cor}[prop]{Corollary}
\newtheorem{theo}[prop]{Theorem}

\newtheorem{ex}[prop]{Example}

\newenvironment{rem}{\refstepcounter{prop}
{\bf{Remark }}\theprop :}{}

\newcommand{\eqnsection}{
   \renewcommand{\theequation}{\thesection.\arabic{equation}}
   \makeatletter
   \csname @addtoreset\endcsname{equation}{section}
   \makeatother}

\def \BE{\hbox{ I\hskip -2pt E}}

\def \PE{\hbox{I\hskip -2pt P}}

\def \N{{\mathbb  N}}
\def\C{{\mathbb C}}

\def \R{{\mathbb R}}

\def\nn{\noindent}
\def\part{\partial}
\def\ra{\rightarrow}

\def\lbc{\lbrace}
\def\rbc{\rbrace}

\def\tilde{\widetilde}

\def\Aa{{\mathcal A}}

\def\prf{{\bf Proof.}}

\def\tr{{\hbox{\rm tr}}}
\def\e{\epsilon}
\def\d{\delta}

\def \E{\hbox{ I\hskip -2pt E}}

\def \PE{\hbox{I\hskip -2pt P}}
\def \NN{\hbox{\it I\hskip -2pt N}}
\def \C{{\mathbb C}}
\def \R {{\mathbb R}}

\def \BE{\hbox{ I\hskip -2pt E}}

\def \PE{\hbox{I\hskip -2pt P}}

\def \N{{\mathbb  N}}
\def\C{{\mathbb C}}

\def \R{{\mathbb R}}

\def\nn{\noindent}
\def\part{\partial}
\def\ra{\rightarrow}

\def\lbc{\lbrace}
\def\rbc{\rbrace}

\def\tilde{\widetilde}

\def\Aa{{\mathcal A}}

\def\Ca{{\mathcal C}}
\def\Ea{{\mathcal E}}

\def\prf{{\bf Proof.}}

\def\tr{{\hbox{\rm tr}}}
\def\e{\epsilon}
\def\tr{{\mbox{tr}}}

\def\NC{{{\mbox{NCP}}}}
\def\cro{{{\mbox{cr}}}}

\def\nn{\noindent}

\def\em{{ empirical measure }}

\def\em{{ empirical measure }}

\def\a{\alpha}

\def\d{\delta}

\def\l{\lambda}

\def\s{\sigma}

\def\e{\epsilon}

\def\arg{\mbox{arg}}

\def\D{\Delta}

\def\nn{\noindent}
\def\part{\partial}

\def\ts{\times}

\def\ra{\rightarrow}

\def\lbc{\lbrace}
\def\rbc{\rbrace}

\def\tilde{\widetilde}

\def\m1{{{\mathcal M}_1^+(\R)}}

\def\prf{{\bf Proof.}}

\newcommand{\beaa}{\begin{eqnarray*}}
\newcommand{\eeaa}{\end{eqnarray*}}
\newcommand{\bea}{\begin{eqnarray}}
\newcommand{\eea}{\end{eqnarray}}
\newcommand{\be}{\begin{equation}}
\newcommand{\ee}{\end{equation}}

\def\11{{\hbox{1\kern-.2em\hbox{I}}}}
\def\pin0{{\hat\pi^N_0}}

\begin{document}

\title{Long time behavior of the
solutions to \\
non-linear Kraichnan equations}
\author{Alice Guionnet \thanks{\'Ecole Normale Sup\'erieure de Lyon,
Unit\'e de Math\'ematiques pures et appliqu\'ees,
UMR 5669,
46 All\'ee d'Italie, 
69364 Lyon Cedex 07, France},\ 
and Christian Mazza\thanks{Section de Math\'ematiques, 2-4 Rue du Li\`evre, CP 240, CH-1211 Gen\`eve 24, Suisse}}
\date{
}

\maketitle

\centerline{\bf Abstract}
We consider the solution of 
a nonlinear Kraichnan equation
$$\partial_s H(s,t)=\int_t^s H(s,u)H(u,t) k(s,u) {\rm d}u,\quad s\ge t$$
with a covariance kernel
$k$ and boundary condition
$H(t,t)=1$. We study
the  long time behaviour of $H$
as the time paramters $t,s$ go to infinity, 
according
to the
asymptotic behaviour
of $k$. This question appears
in various
subjects since it 
is related with the analysis 
of the asymptotic behaviour
of the trace 
 of 
non-commutative 
processes satisfying  a linear differential equation,
but also naturally shows up 
in the study  of the
so-called response function and aging
properties
of the dynamics of some
disordered spin systems.

{\vspace{1cm}}
\nn

\nn
{\underline{\bf  Mathematics Subject of Classification  :}}
 82B44, 46L54, 45G10.

\section{Introduction\label{s.introduction}}


In this paper, we
shall consider 
the long time behaviour
of the solution 
of a nonlinear Kraichnan equation
\begin{equation}\label{kr}
\partial_s H(s,t)=\int_t^s H(s,u)H(u,t) k(s,u) {\rm d}u
\end{equation}
with a covariance kernel
$k$ and boundary condition
$H(t,t)=1$.
Such an equation already appeared
in the work of Kraichnan  \cite{Kr}
as a first term in a perturbative method
to analyze quantum-mechanical, turbulence
or disordered  problems. Shortly afterwards,
Frisch and Bourret \cite{FB} have shown
that these equations naturally appeared 
when one considers parastochastic
equations, which are related 
with differential equations
for non-commutative
processes (so-called master equations)
and large random matrices.
This relation was later studied also 
by Neu and Speicher \cite{NS}. Let us briefly describe
it.

\smallskip

Let $(L_t)_{t\ge 0}$ be
a process in a von Neumann algebra $\Aa$
equipped with a tracial state $\phi$.
We assume that $L$ is a centered semicircular
process with covariance kernel
$k$, usually constructed on the
full Fock space (see e.g. \cite{Voi}).
In a more intuitive 
 way, $L$ can be constructed as the limit
 of self-adjoint large random matrices
$(L^N_t)_{t\ge 0}$ with  entries
$\{(L^N_t)_{ij}, 1\le i\le j\le N\}$
which are independent 
 Gaussian processes
with covariance $N^{-1} k$. This limit has to be understood
in the weak sense that for any integer number $n$,
any times $(t_1, t_2,\cdots t_n)\in
(\R^+)^n$,
$$\lim_{N\ra\infty} {1\over N}\tr\left(L^N_{t_1}L^N_{t_2}
\cdots L^N_{t_n}\right)=\phi\left(
L_{t_1}L_{t_2}
\cdots L_{t_n}\right),$$
where $\tr$ denotes the unnormalized  trace of matrices.
\medskip

\nn
In this paper, we consider
the family of operators
${\bf X}_{s,t}$ satisfying the linear differential equation
$$\partial_s {\bf X}_{s,t} =L_s
 {\bf X}_{s,t},\ s>t,$$
with boundary data ${\bf X}_{t,t}=1$, on the full Fock space.
Then, it was shown in \cite{FB}
that $H(s,t)=\phi({\bf X}_{s,t})$
satisfies Kraichnan's equation (\ref{kr}) (see also
section \ref{s.s.formal} for details).
We  study the asymptotic behaviour of
$H(s,t)=\phi({\bf X}_{s,t})$ as $s$ and $t$  go to infinity. Notice
that $\phi({\bf X}_{s,t})$ describes the large $N$ limit of
the mean normalized trace $N^{-1}\tr({\bf X}_{s,t}^N)$ of
the solution ${\bf X}_{s,t}^N$ of 
the random     linear differential equation
$\partial_s {\bf X}_{s,t}^N
 =L_s^N
 {\bf X}_{s,t}^N,\ s>t$ such that ${\bf X}_{t,t}^N=I$.

\medskip

Such a question would of course be trivial in 
the classical setting where $L$ would just be 
a real-valued  Gaussian process. Indeed, in this
case, ${\bf X}_{t,s}=e^{\int_t^s L_u {\rm d}u}$ and one easily studies
the asymptotics of $\phi({\bf X}_{s,t})=\E[{\bf X}_{t,s}]$ thanks to the formula
\begin{equation}\label{cl}
\E[{\bf X}_{t,s}]=e^{{1\over 2}\int_t^s\int_t^s k(u,v) {\rm d}u {\rm d}v}.
\end{equation}

\medskip

It appears  to
be actually quite a difficult 
question in the non-commutative 
setting. Eventhough it is a rather natural question to
address concerning differential equations in 
free probability, our first
motivation
came actually from standard 
statistical mechanics, namely from the
study of the aging
properties of p-spherical spin glasses.
Indeed, 
consider a spin glass with Hamiltonian
$$H_J(x)=\sum_{p=1}^M \frac{a_p}{p!}\sum_{1\le i_1\cdots i_p\le 
N}J_{i_1\cdots i_p}x^{i_1}\cdots x^{i_p},$$
$x=(x^i)_{1\le i\le N}$, $x^i\in {\R}$,
with independent centered gaussian variables $J_{i_1\cdots i_p}$. 
The Langevin dynamics for this model 
with a smooth spherical constraint 
are given by  the stochastic differential system

$${\rm d}x_t = -f'(\vert\vert x_t\vert\vert^2/N)x_t {\rm d}t - \nabla 
H_J(x_t){\rm d}t +{\rm d}B_t,$$
where $B_t$ is a $N$-dimensional Brownian motion and $f$ is a convex
 function. Let 
$$C_N(s,t)=\frac{1}{N}\sum_{i=1}^N x_s^i x_t^i,\ \ {\mathcal 
X}_N(s,t)=\frac{1}{N}\sum_{i=1}^N x_s^i B_t^i.$$
 It was shown in \cite{BDG2}
that the couple $(C_N,\chi_N)$
converges almost surely  towards functions $(C,\chi)$.
If we set

$$R(s,t)=\partial_s {\mathcal X}(s,t)$$
be the so-called response function of the system,
then (see \cite{BDG2}), $(C,G)$  
 satisfies the following
integro-differential equations given, for $t\le s$,  by
$$\partial_s R(s,t)=-f'(C(s,s))R(s,t)+\int_t^s
R(s,u)R(u,t)\nu''(C(s,u)){\rm 
d}u,$$
\begin{eqnarray*}
\partial_s C(s,t)=-f'(C(s,s))C(s, 
t)&+&\int_0^t R(t,u)\nu'(C(s,u)){\rm d}u\\
&+&\int_0^s R(s,u)C(t,u)\nu''(C(s,u)){\rm d}u,
\end{eqnarray*}
where  the function $\nu$ is given by
$$\nu(x)=\sum_{p=1}^M \frac{a_p^2}{p!} x^{p}.$$
Here, the boundary conditions are given by $R(t,t)\equiv 1$ and
$C(0,0)$ (which is known).
Similar type of equations
have been derived in various contexts such as
the dynamics of long-range superconducting
networks \cite{Cha}
or for other dynamical models \cite{Fra}.

\medskip

The main question which arises
in physics is to understand the
long time behaviour of $C$, which measures
the long time memory of the system 
and aging (see \cite{BDG} for a detailed study of the
easiest case $\nu(x)={c\over 2}x^2$). 
  \cite{CuKu}
derived the same set of 
Schwinger-Dyson equations
(for the hard sphere
model where $f$ is a function of the time variable only,
chosen  so that
$C(t,t)\equiv 1$ for all $t\in\R^+$) and 
 proposed heuristics concerning
the asymptotic behavior of the solutions when
$t$ and $s$ are large
(see also \cite{Cu}). However, even 
on a non-rigorous ground, no complete 
description of these asymptotics could
be given, but rather scenarios 
about their form could be validated or dismissed. 
The idea is indeed to assume a form for the asymptotics
of the couple $(C,R)$ in different 
regimes of the time parameters $(t,s)$ ;
for instance, one can imagine that
on $D_{\mbox{FDT}}:=\{ s\ge t : s-t\gg t\}$ 
(corresponding to the 
so-called FDT regime), the solutions are approximately 
stationary
$$C(s,t)\equiv C_{\mbox{FDT}} (s-t),\quad R(s,t)\equiv 
R_{\mbox{FDT}} (s-t),$$
with a standard choice of the form $C_{\mbox{FDT}}(x)=Ae^{-ax},
R_{\mbox{FDT}} (x)=Be^{-ax}$,
but on $D_{\mbox{AGING}} =\{ s\ge t : t/s\in (0,1)\}$
 (corresponding to an aging
regime), one can expect
$$C(s,t)\equiv C_{\mbox{AGING}} ({t\over s}
),\quad R(s,t)\equiv {1\over s} R_{\mbox{AGING}} ({t\over s}
)$$
with a standard guess $C_{\mbox{AGING}}(x)=Ax^p,
R_{\mbox{AGING}} (x)=Bx^q$ for some exponents 
$p,q$ to be determined. 
Then, one checks whether this scenario
is consistent with 
the above integro-differential system.
However, the form of the intermediate phases
between these different domains 
is hard to predict and actually determines the exponents
(such as $p,q$)
since they give the slope of the different curves at
their boundary (in fact, the integrals
in our system of equations
keep track of
all the past of the trajectories, including
these intermediate phases). This is
the so-called matching problem (between different phases). 
 Hence, such a strategy can not, 
without an intuition
about these intermediate phases, predict
completely the solution.

\medskip
In this paper, we shall study only the equation
for the response function being given the
covariance $C$, which asymptotics shall take the 
forms that we expect to encounter for
the solution of the full system.
 Hopefully,
this will enable us 
to perform later a bootstrap argument 
to solve our original problem 
concerning the asymptotic behaviour
of the covariances of p-spherical systems for instance.
At least, we hope it will shade some light
on the behaviour of the response function.

\medskip

If we set 
$$H(s,t)=\exp(\int_t^s f'(C(u,u)){\rm d}u)R(s,t),$$
it is easy to check that $H$ solves the equation (\ref{kr})
with $k(s,u)=\nu''(C(s,u))$.  Hence, studying the asymptotics 
of the  response function
being given 
the covariance $k$ is equivalent to
study of the long time behaviour of the solution of (\ref{kr}).
As a remark, we want to point
out that we  have no deep insight
 why the response function should be related 
with a non-commutative process
; we only realize that both evolutions
 are given by
the
same integro-differential equation. It is not clear if
the full system could have such an operator interpretation.

\medskip

Amazingly, our motivation brings
us back to Kraichnan \cite{Kr}
who considered the equations (\ref{kr}) 
when trying to analyze
the Schr\"odinger equation of a particle in a random 
potential; his method 
provides such equations for correlation 
functions and averaged Green's
functions ! This coincidence 
might lie in the fact that 
spherical models are well
suited for his expansion
method, but we shall
not  study this
question here.

\medskip

Throughout this article we will assume that

\noindent {\bf Hypothesis:}
We shall assume that $k$ is non negative
and uniformly
bounded, 
i.e that
\begin{eqnarray}\label{hypo}
k(t,s)\ge 0,\quad t,s\in\R^+\quad, \sup_{s,t\in\R^+}
k(t,s)=\sup_{s\in\R^+}
k(s,s)<\infty.
\end{eqnarray}
This hypothesis should be fulfilled by the covariances of the
p-spins models.

Let us now state our main results.
Our most precise estimates 
are obtained in the
case where the covariance is stationary,
in which case the Laplace transform of the solution
$H$ to (\ref{kr}) is given by rather nice formulas
so that we can use complex analysis to
estimate $H$.

We have then the following
dichotomy depending  whether
$k$ converges to
zero or not as time
goes to infinity ;

\begin{theo}\label{main}
Suppose that $k$ is stationary, that is 
$k(t,s)=k(s-t)$. Then
1)
Assume that the kernel $k$ is
 such that there exists $a>1$ and $C<\infty$
so that
$$0\le k(u)\le {C\over (1+u)^{a}}.$$
Then, there exists $\l_c(H)>0$ such that
$$
\exp(-\l_c(H) t)H(t)\sim \frac{1}{2A},\ \ t\to +\infty,$$ 
where
$$A=\frac{{\rm d}}{{\rm d}\l}(\l-\hat{Hk}(\l))\vert_{\l=\l_c(H)}>0.$$

2)Assume that
$k(u)=c_2+c_1 k_1(u)$, for some
positive constants $c_1$, $c_2$ and $|k_1(u)|\le C(1+u)^{-1}$
for some finite $C$.
Then, there exists $\l_c(H)>0$ such that
$$e^{-\l_c(H) t} H(t)\sim A t^{-3/2},\ \ t\to+\infty,$$
for some positive constant $A$.
\end{theo}
This theorem is proved in Theorems \ref{vancov}
and \ref{algdec}. We can not in general compute the 
Lyapounov exponent $\l_c(H)$ except
in the case where $k(u)=ce^{-\d u}$.
In this case, which we study
in details
in section \ref{expvan}, $\l_c(H)$ 
appears to be  the 
 smallest zero of
a Bessel function. Interestingly,
this special case has already 
been studied for combinatorial  reasons
in \cite{Brak} ; the Laplace 
transform of $H$ can in fact be interpreted as
 the generating function
of random staircase polygons.

We also consider the
general case and prove that

\begin{theo}\label{flatintro}
Assume that
$$k(s,t)= k_1(s-t)+h(s,t)$$
with a flat function $h$ 
such that
there exists a positive constant $C$
and for $T,M>0$ a function $\d(T,M)$ such that
$\d(T,M)\ra 0$ for all $M$ when $T$ goes to
infinity
so that

\begin{equation}\label{hypint}
\sup_{t\ge T}\sup_{|s-t|\le M}| h(s,t)-C|\le
\d(M,T),\quad \sup_{s\ge t\ge T}|h(s,t)|\le C.
\end{equation}
Then, we claim that if we denote $H_k$
the solution of (\ref{kr}), 
regardless of
the way $t\ge T$ goes to
infinity 
$$\lim_{t\ra\infty}\lim_{s-t\ra\infty}{1\over s-t}\log
H_k(s,t)=\lim_{t\ra\infty}\lim_{s-t\ra\infty}{1\over s-t}\log
H_{C+k_1} (s,t)=\l_c(H_{C+k_1}).$$
\end{theo}
The study of the second order 
correction in this
general case is of course much more complicated
and demand much more precise hypothesis 
concerning the function $h$.
We shall come back to
this issue in a forthcoming research.

\bigskip

The plan of the article is as follows ;
we first show that
 $(\phi({\bf X}_{s,t}), s\le t)$
can be described as  the unique solution
of
(\ref{kr}) in section \ref{s.s.formal},
we then study the robustness
of the first order
asymptotics of $H$ when $k$ varies in section \ref{weaks}.
In section \ref{stationary},
we consider the case where $k$ is stationary.
The general case is considered in section \ref{general}.

\section{Description of $(\phi({\bf X}_{s,t}), s\le t)$ as
the unique solution of the integro-differential
equation } \label{s.s.formal}

Let $L$ be a semicircular 
process with covariance $k$
and ${\bf X}_{s,t}$ satisfying the linear differential equation
$$\partial_s {\bf X}_{s,t} =L_s {\bf X}_{s,t},$$
with boundary data ${\bf X}_{t,t}=1$.
We first remind the reader why $(\phi({\bf X}_{s,t}),s\ge t)$
satisfies (\ref{kr}) and then that
it is actually uniquely described by this property.

\subsection{$(\phi({\bf X}_{s,t}),s\ge t)$
satisfies (\ref{kr})}

Let us recall 
that a semi-circular variable $L_t$
with covariance $k(t,t)$ is uniformly bounded
(for the operator norm) by $2k(t,t)^{1\over 2}$.
Therefore, we can use 
Picard  formula, which
can serve as a basic definition of ${\bf X}_{s,t}$,
to write ${\bf X}_{s,t}$ in the form
\begin{equation}\label{picard}
{\bf X}_{s,t}=\sum_{n\ge 0}
\int_{t\le t_1\cdots \le t_{n}\le s}
L_{t_1}\cdots L_{t_{n}}
{\rm d}t_1\cdots {\rm d}t_n
\end{equation}
where the serie converges
uniformly with respect to
the operator norm and uniformly on  any times $s,t$ in a compact
interval since $k$ is uniformly bounded.

In a combinatorial way,
the fact that $L$ is
a semicircular
process means that the free cumulants
of $L$ are null except the second one,
which is given by $k$. 
The
analogue of Wick formula for
such processes is given
by
\begin{equation}\label{wick}
\phi(L_{t_1}L_{t_2}...L_{t_{2n}})
=\sum_{\sigma\in \NC_n}\prod_{i\in
\cro(\s)} k(t_i,t_{\s(i)}),
\end{equation}
where $\NC_n$ denotes the
set of involutions of $\{1,\cdots,2n\}$
without fixed points and without crossings
and where  $\cro(\sigma)$
is defined to be the set of indices $1\le i\le 2n$
such that $i<\sigma(i)$.  $\sigma\in \NC_n$
when the situation $i<j<\sigma(i)<\sigma(j)$ does not
occur.
  Therefore,  $B(t_1,\cdots,t_n):=\phi(L_{t_1}\cdots L_{t_{n}})$
is null when $n$ is odd and otherwise
satisfies the recursion formula
$$
B(t_1,\cdots, t_{2n})=\sum_{p=1}^{2n}
k(t_1,t_p)B(t_2....t_{p-1})B(t_{p+1}...t_{2n}).
$$
As a consequence of (\ref{picard}), 
we get

\begin{lem}\label{formalsolution}
\begin{equation}\label{generalformula}
H(s,t)=\phi({\bf X}_{s,t})=\sum_{n\ge 0}
\int_{t\le t_1\cdots \le t_{2n}\le s}
\sum_{\sigma\in \NC_n}\prod_{i\in
\cro(\s)} k(t_i,t_{\s(i)}){\rm d}t_1\cdots {\rm d}t_{2n}
\end{equation}
 solves (\ref{kr}) and satisfies
\begin{equation}\label{bound1}
H(s,t)\le \exp(2\int_t^s k(u,u)^{1/2}{\rm d}u).
\end{equation}
\end{lem}
\prf 
Indeed, by definition, $\partial_s H(s,t)$ is given by
$$
\sum_{n\ge 0}\sum_{i:\s(i)=2n}
\int_{t\le t_1\cdots \le t_{2n-1}\le s} k(t_i,s)
\sum_{\sigma\in \NC_n\backslash \{2i-1,2n\}}\prod_{i\in
\cro(\s)} k(t_i,t_{\s(i)})\prod_{i=1}^{2n-1}{\rm d}t_i$$
$$
=\sum_{n\ge 0}\sum_{i=1}^{n}
\int_{t\le t_{2i-1}\le s} k(t_i,s)$$
$$\left( 
\int_{t\le t_1\cdots \le t_{2i-2}\le t_{2i-1}}
\sum_{\sigma\in \NC_{i-1}} 
\prod_{j\in
\cro(\s)} k(t_j,t_{\s(j)}) {\rm d}t_1\cdots {\rm d}t_{2i-2}\right)$$
$$
\ts
\left( 
\int_{t_{2i-1} \le t_1\cdots \le t_{2(n-i-1)}\le s}
\sum_{\sigma\in \NC_{n-i-1}} 
\prod_{j\in
\cro(\s)} k(t_j,t_{\s(j)}) {\rm d}t_1\cdots {\rm d}t_{2(n-i)-2}\right) {\rm d}t_{2i-1}
$$
$$= \int_t^s k(u,s)H(u,t)  H(s,u) {\rm d}u,$$
where in the second line we noticed 
that $\{i:\s(i)=2n\}=\{1,3,\cdots, 2n-1\}$
since $\s\in \NC_n$ and we obtained the last one by 
summing over the indices $1\le i\le n\le\infty$.
The only point which remains to prove is (\ref{bound1}) :
\begin{eqnarray*}
\phi({\bf X}_{s,t})&=& \sum_{n\ge 0}
\sum_{\s\in\NC_n}\int_{t\le t_1\cdots\le t_{2n}\le s}
\prod_{i\in \cro(\s)}
k_2(t_i,t_{\s(i)})
\prod {\rm d}t_i\\
&\le&\sum_{n\ge 0}\sum_{\s\in\NC_n}
{1\over 2n !}
\left( \int_t^s
   k(u,u)^{1\over 2} {\rm d}u\right)^{2n}\\
&=& \sum_{n\ge 0} C_n
{1\over 2n!}\left(\int_t^s
   k(u,u)^{1\over 2} {\rm d}u\right)^{2n}
\\
&=&E(e^{\int_t^s k(u,u)^{1\over 2} {\rm d}u S}),\\
\end{eqnarray*}
where we used Cauchy-Schwarz inequality, where $C_n$
denotes the Catalan number of order $n$
(i.e the number of partitions in $\NC_n$),
and $S$ a standard semi-circular random variable
(i.e. a random variable with
law $\s(dx)=(2\pi)^{-1}\sqrt{4-x^2} dx$) 
which is well known to satisfy $\E[S^{2n}]=C_n$.
Using the fact that $S$ is bounded
by $2$ uniformly, we obtain (\ref{bound1}).

\hfill\endproof

\begin{rem}\label{bound3}
When $k$
is uniformly bounded, say by $C$, 
we deduce 

$$ \phi({\bf X}_{s,t})\le e^{2\sqrt{C}(s-t)}.$$
Moreover, a 
look at the previous proof shows that,
if $k(t,s)\equiv C$,
$$\lim_{s-t\ra\infty}{1\over s-t}\log \phi({\bf X}_{s,t})
=2\sqrt{C}.$$
\end{rem}

\bigskip


\subsection{Uniqueness of the solution of (\ref{kr})}
Set
$$\Ea_M=\{ f\in\Ca_b(\R^+\ts\R^+):\quad |f(s,t)|\le Me^{M|t-s|}\}
$$

\begin{theo}\label{unique}
Assume that $k$ is uniformly bounded
on $\R^+\ts\R^+$ by some constant $C$. Then, for any $M\in\R^+$, there exists
at most one solution to (\ref{kr}) in
$\Ea_M$. Moreover if $M\ge 2\sqrt{C}$,
the solution is given by (\ref{generalformula}).
\end{theo}
\prf
  Consider two solutions in $\Ea_M$, $H$ and
$\tilde H$ and set
$$\D(s,t)=e^{-M|s-t|} |H(s,t)-\tilde H(s,t)|.$$
Then
\begin{eqnarray*}
\D(T,t)&\le& CM\int_{t\le u\le s\le T}
[\D(s,u)+\D(u,t)] e^{M(s-T)} {\rm d}u{\rm d}s\\
&\le &  CM\int_{t\le u\le s\le T} \D(s,u) e^{M(s-T)} {\rm d}u{\rm d}s
+C\int_{t\le u\le T}\D(u,t) {\rm d}u
\\
\end{eqnarray*}
Using Gronwall's lemma ($\D$ is bounded by hypothesis)
   to get rid of  the last term in the above right hand
side, we deduce
\begin{eqnarray}
\D(T,t)&\le& C(e^{C(T-t)}(T-t+1) 
\int_{t\le u\le s\le T}\D(s,u) e^{-M(T-s)}{\rm d}u{\rm d}s
\label{b1}\\
\nonumber
\end{eqnarray}
Consequently, $F(T)=\sup_{0\le t\le T} \D(T,t)$
satisfies
$$F(T)\le C (T+1)e^{CT} \int_0^T  
 F(s)  {\rm d}s $$ 
is null by Gronwall's lemma,
resulting with $\D$ null by (\ref{b1}).
The last part of the statement is a consequence
of Remark \ref{bound3}.

\hfill\endproof

\begin{example}\label{Decoupled}
Suppose that
$k(t,s)=h(s)h(t)$, for some function $f$.
In such a case, we simply take $L_t=h(t) S$
with a given semicircular variable
$S$. Then,  the
solution is the 'classical' one
$${\bf X}_{s,t}=e^{ S\int_t^s h(u) {\rm d}u},$$
and therefore,
\begin{eqnarray*}
H(s,t)&=&{\rm E}(e^{ S\int_t^s h(u) {\rm d}u})\\
&=&(2\pi)^{-1}\int e^{ \int_t^s h(u){\rm d}u x} \sqrt{4-x^2}dx.
\end{eqnarray*}
Consequently, when $\int_t^s h(u){\rm d}u$
goes to infinity,
$$H(s,t)\approx (\int_t^s h(u){\rm d}u)^{-3/2}\exp\{2 \int_t^s h(u){\rm d}u\}.$$
This can be compared with 
the classical setting where
(see (\ref{cl}))
$$\BE[{\bf X}_{s,t}]=e^{{1\over 2}(\int_t^s h(u) {\rm d}u)^2}.$$
\end{example}


\section{Weak continuity statements}\label{weaks}


In this section, we shall investigate 
the robustness 
of the asymptotic behaviour
of $H$ when $k$ varies. 
Let us first note that by (\ref{generalformula}),
it is clear that
\begin{pr}\label{increas}
For any $t_0\ge 0$, any covariance kernels $k_1,k_2$
such that $0\le k_1(s,t)\le k_2(s,t)$
for all $s\ge t\ge t_0$,
$$H_{k_1}(s,t)\le H_{k_2}(s,t),\quad \forall s\ge t\ge t_0$$
\end{pr}
\begin{prop}\label{weak}
Let $k_1(s,t)$ and $k_2(s,t)$ be two covariance
functions such that
  for any $\e>0$, there is $t_\e<\infty$
such that for $s\ge t\ge t_\e$
$$(1-\e) k_2(s,t)\le k_1(s,t)\le (1+\e)k_2(s,t).$$
Then, denoting by $H_k$  the
solution of (\ref{kr}) with kernel $k$, we have,
uniformly for  any $t>t_\e$
$$\liminf_{s-t\ra\infty}
{1\over s-t}\log H_{k_2}(s,t)+2Ce\log(1-\e)
\le
\liminf_{s-t\ra\infty}
{1\over s-t}\log H_{k_1}(s,t)\le$$
$$\le\limsup_{s-t\ra\infty}
{1\over s-t}\log H_{k_1}(s,t)\le 
 \limsup_{s-t\ra\infty}
{1\over s-t}\log H_{k_2}(s,t)+2Ce\log(1+\e),$$
where $C=\sup_{s\in\R^+}k(s,s)^{1/2}$.
\end{prop}
\prf
 From  (\ref{generalformula}), we know that
$$H_k(s,t)=\sum_{n\ge 0} B_n^k (s,t)$$
with
$$B_n^k(s,t)=\int_{t\le t_1\le\cdots\le t_{2n}\le s}
\phi(L_{t_1}\cdots L_{t_{2n}}) \prod {\rm d}t_i.$$
Hence, if $A=\sum_{n\ge 1} n^{-2}$,
\begin{equation}\label{bound2}
\max_{n\ge 0}(B_n^k(s,t))\le H_k(s,t)\le A
\max_{n\ge 0} ((n+1)^2 B_n^k(s,t)).
\end{equation}
We already noticed that
$$B_n^k(s,t)\le {(\int_t^s k(u,u)^{1\over 2} {\rm d}u)^{2n}\over (2n)!}C_n
\le {(2C(s-t))^{2n}\over(2n)!}$$
where $C$ is a bound on $k^{1\over 2}$
and where we used $C_n\le 4^n$.
As a consequence, using Stirling formula, for any $B>2C$,
and $s-t$ large enough,
$$\max_{n\ge B(s-t)}(n+1)^2 B_n^k(s,t)\le (B(s-t)+1)^2
\left({e C\over B}\right)^{2
B(s-t)}.$$
Therefore, 
fixing any $B>eC$, say $B=2eC$,
we see that there exists $M<\infty$ such that 

$$\sup_{s-t>M\atop t\in\R^+}
\max_{n\ge 2eC(s-t)}(n+1)^2 B_n^k(s,t) \le 1.$$
Consequently, $\sup_{s-t>M\atop t\in\R^+}
\max_{n\ge 0} ((n+1)^2 B_n^k(s,t))$ is given by
$$
\sup_{s-t>M\atop t\in\R^+}
\max\left\lbc \max_{n\le 2eC(s-t)} ((n+1)^2 B_n^k(s,t)),
\max_{n\ge 2eC(s-t)} ((n+1)^2 B_n^k(s,t))\right\rbc$$
$$
\le \max\left\lbc
\sup_{s-t>M\atop t\in\R^+}
\max_{n\le 2eC(s-t)} ((n+1)^2 B_n^k(s,t)), 1\right\rbc.$$
But,
 by definition $B_0^k(s,t)\equiv 1$
so that in fact
$$\sup_{s-t>M\atop t\in\R^+}
\max_{n\ge 0} ((n+1)^2 B_n^k(s,t))\ge 1$$
and therefore for any $s-t>M$, any $t\in\R^+$,
$$
\max_{n\ge 0} ((n+1)^2 B_n^k(s,t))=
\max_{n\le 2eC(s-t)} ((n+1)^2 B_n^k(s,t)).$$
As a consequence, we deduce from
(\ref{bound2}) that for any $s-t>M$, any $t\in\R^+$,

\begin{equation}\label{kl}
\max_{n\ge 0}(B_n^k(s,t))\le H_k(s,t)\le A (2eC(s-t))^2
\max_{n\le 2eC(s-t)}  B_n^k(s,t).
\end{equation}
We thus deduce that,
regardless of the way $t$ goes to infinity
(or not),
\begin{equation}\label{egal1}
\limsup_{s-t\ra\infty}
{1\over s-t}\log H_k(s,t)= \limsup_{s-t\ra\infty}
{1\over s-t}\log \max_{n\le 2eC(s-t)}  B_n^k(s,t),
\end{equation}
and
$$
\liminf_{s-t\ra\infty}
{1\over s-t}\log H_k(s,t)= \liminf_{s-t\ra\infty}
{1\over s-t}\log \max_{n\le 2eC(s-t)}  B_n^k(s,t).
$$
Now, for $t\ge t_\e$, our hypothesis implies 
 $B_n^{k_1}(s,t)\le (1+\e)^n B_n^{k_2}(s,t)$, which 
results with
$$\limsup_{s-t\ra\infty}{1\over s-t}\log H_{k_1}(s,t)\le
2Ce \log(1+\e) +\limsup_{s-t\ra\infty}{1\over s-t}\log H_{k_2}(s,t)$$
The same arguments apply for the lower bound.

\hfill \endproof

We now show a slightly stronger result 
giving the first order asymptotics
of $H$ for slowly decaying covariances

\begin{cor}\label{slow}
Suppose that $k$ is a covariance
such that for all $\e>0$,
there exists $t_\e<\infty$
such that for all $s\ge t\ge t_\e$,
\begin{equation}\label{as}
Ce^{-\e(s-t)}\le k(t,s)\le C ,
\end{equation}
for some positive constant $C>0$.
Then
$$\lim_{t\ra\infty}\lim_{s-t\ra\infty}{1\over s-t}\log H_k(s,t)
=2\sqrt{C},$$
corresponding to the limit where $k(s,t)=C$
(see Example \ref{Decoupled}).
Moreover, if we assume additionally that
$k(s,t)$ is  decreasing in $s$ 
and increasing in $t$ on $s\ge t$, we also have 
that for any $\d>0$ there exists $M_\d<\infty$
so that for $(s,t)$ such that $(s-t)\sqrt{k(s,t)}\ge M_\d$,
$${1\over s-t}\log H_k(s,t)\ge (2-\d)\sqrt{k(s,t)}.$$
Thus, (\ref{as}) implies

$$
\lim_{s-t\ra\infty}\lim_{t\ra\infty}{1\over s-t}\log H_k(s,t)
=2\sqrt{C}.$$
More precisely, we have in general
$$ \lim_{s-t\ra\infty}\sup_{k(s,t)\ra 1} {1\over s-t}\log H_k(s,t)
=\lim_{s-t\ra\infty}\inf_{k(s,t)\ra 1} {1\over s-t}\log
H_k(s,t)=2\sqrt{C}
.$$
\end{cor}
\begin{rem}\label{ratio}
Take $k(s,t)=C(t/s)^\a$ for $s\ge t$ 
for some $\a>0$.
Then, it is easy to check
that for any $t\ge t_\e={\a\over\e}\sup_{v>0}{1\over v}\log(1+v)$,
any $s\ge t$,
$$ Ce^{-\e(s-t)}\le  k(s,t)\le C$$
so that the conclusions of Corollary \ref{slow} apply.
Note that this corollary only concerns the cases where 
$s-t$ and $t$ {\bf go to infinity independently}
or at most in such a way that $(t/s)$ goes to
one. In such regimes,
$k$ converges either to zero (when $s-t$ goes to
infinity first) or one (when $t$ goes to infinity first).
We shall consider in section \ref{general} 
the case where (\ref{as}) is generalized 
to the case where $C$ is not constant
but a stationary
function. 
\end{rem}

\nn
\prf
By Property \ref{increas},
we  see that for $s\ge t\ge t_\e$,
$$H_{C e^{-\e(s-t)}}(s,t)\le H_k(s,t)\le H_C(s,t).$$
Therefore,
$$\liminf_{\e\ra\infty}\liminf_{s-t\ra\infty}{1\over s-t}\log
H_{C e^{-\e(s-t)}}(s,t)\le 
\liminf_{t\ra\infty}\liminf_{s-t\ra\infty}{1\over s-t}\log 
 H_k(s,t)\le$$
$$
\le \limsup_{t\ra\infty}\limsup_{s-t\ra\infty} 
{1\over s-t}\log H_k(s,t)\le \limsup_{t\ra\infty}
\limsup_{s-t\ra\infty} 
{1\over s-t}\log H_C(s,t)=2\sqrt{C}$$
where the first equality comes
from the observation
that 
$H_{C e^{-\e(s-t)}}(s,t)=H_{C e^{-\e(s-t)}}(s-t)$
so that taking $t$ large only results in taking
$\e$ as small as wished and the
last equality comes from Remark \ref{bound3}.
We shall
see in Proposition \ref{Bessel}
that for any $\e>0$,
$$\liminf_{s-t\ra\infty}{1\over s-t}\log
H_{C e^{-\e(s-t)}}(s,t)=\l_c(\e)$$
and further that $
\l_c(\e)$ converges towards $2\sqrt{C}$
as $\e$ goes to zero, finishing the proof of our first result.

\medskip

There is  an easy
argument to prove directly
the second  statement ;
assume that $k(s,t)$
decreases in $s$
for all $t$ and increases in $t$  so that $k(u,v)\ge k(t,s)$
for all $t\le u\le v\le s$.
Then, by (\ref{generalformula}), we find that
\begin{eqnarray*}
H(s,t)&\ge&\sum_{n\ge 0}\sum_{\s\in\NC_n}
k(s,t)^n {(s-t)^{2n}\over (2n)!}\\
&=&\sum_{n\ge 0} C_n 
k(s,t)^n {(s-t)^{2n}\over (2n)!}=\E[ e^{S(s-t)\sqrt{k(s,t)}}]
\end{eqnarray*}
with a semicircular variable $S$. Thus, for any $\d>0$,

$$H(s,t) \ge \PE(S>2-\d)e^{(2-\d)(s-t)\sqrt{k(s,t)}}$$
yielding the estimate since $\PE(S>2-\d)>0$ for any $\d>0$.
As a consequence, we trivially get the last point of
the corollary
 since we already have the upper
bound.

\hfill \endproof

\section{Asymptotic behaviour of $H$ for stationary covariances}
\label{stationary}

When $k(t,s)=k(s-t)$, (\ref{generalformula}) yields
$$
H(s,t)=
\sum_{n\ge 0}\sum_{\s\in\NC_n} \int_{0\le t_1
\le\cdots\le t_{2n}\le s-t}  \prod k(t_{\s(i)}-t_i) {\rm d}t_j$$
so that $H(s,t)=H(s-t)$, and
  (\ref{kr}) becomes
\begin{equation}\label{f3}\partial_t H(t)=\int_0^t H(t-u)H(u)k(t-u){\rm d}u
\end{equation}
Consider  the (eventually infinite)
 Laplace transform 
$$\hat H(\l)=\int_0^\infty e^{-\l u}H(u) {\rm d}u.$$
Observe that since $H(u)\ge 1>0$ for
all $u\in\R^+$, $\hat H$ is strictly decreasing.
Moreover,  (\ref{bound1})
(see Remark \ref{bound3}) shows
 that
$\hat H$ is finite for $\lambda$
large enough. The region of
convergence of the Laplace transform is of
the form $(\lambda_c,+\infty)$, for some
critical parameter $\lambda_c$. By assumption, the kernel 
$k$ is non-negative, implying that
$\hat H$ diverges to $+\infty$ on
$(-\infty,\lambda_c)$.
$\hat H$ is analytic on its domain
of convergence and the non-negativity of
$k$ implies that the abscissa of
convergence $\l_c$ is a  singularity
of $\hat H$ (see Theorems 5a and 5b in \cite{Wi}).
Let
$$\l_c(H)=\inf\{\l\in\R : \hat H(\l)<\infty\}.$$
Note that $\l_c(H)<\infty$ by Remark \ref{bound3}
(in fact, $\hat H(+\infty)=0$) and that $\l_c(H)\ge 0$
since $H\ge 1$ so that $\hat H(0)=+\infty$.

Note that since $k$ is uniformly bounded,
$$\l_c(Hk)\le \l_c(H).$$
Moreover, by (\ref{f3})
and using Fubini's theorem for non negative functions,
we find that for any
$\l> \l_c(H)$, 

\begin{eqnarray}
\l \hat H(\l)&=&1+\hat H(\l) \hat{(Hk)}(\l)\label{f5}\\
\nonumber
\end{eqnarray}
We shall now show that depending 
whether  $k$ 
goes  to
zero or not at infinity, the asymptotic
behaviour of $H$  will be rather different.
All the proofs
are based on a refinement of Tauberian
theorems based on   analytic continuations
of the function $\hat H$ 
and the following Lemma 7.2
of \cite{BDG}:
\begin{lem}\label{Inverse}
Suppose that the Laplace transform
$$\hat f(z)=\int_0^\infty e^{-zx}f(x){\rm d}x,$$
of an absolutely  integrable, continuous function
$f(x)$, defined for $\Re(z)>0$, has an analytic continuation
on a domain $S_\theta$ of the form
$$S_\theta=\{z\in\C^*;\ \vert {\rm arg}(z)\vert<(\pi/2)+\theta\},$$
for some
 $\theta\in (0,\pi/2)$, and is such that
$\vert \hat f(z)\vert\to 0$ as $\vert z\vert\to\infty$ in $S_\theta$.
If for some $r>0$ and $A>0$,
$$\limsup_{\vert s\vert\to 0,\ s\in S_\theta}\vert s^r
\hat f(s)-A\vert=0,$$
 then
$$\limsup_{x\to+\infty}\vert x^{1-r}f(x)-\frac{A}{\Gamma(r)}\vert=0.$$
\end{lem}

\subsection{Vanishing covariances }
We start with a
specific example where
we can even precise the value of 
$\l_c(H)$, that is the case 
of exponentially vanishing covariances.
We then tackle the general
case.

\subsubsection{Exponentially vanishing covariances}\label{expvan}

We shall  study precisely 
the asymptotics of $H$ in the case where
 $k(u)=ce^{-\d u}$, for some
$c>0$, $\d>0$. 
We denote here in short $\l_c(\d)=\l_c(H_{ce^{-\d .}})$.
 (\ref{f5}) gives, for $\l> \l_c(\d)$,
\begin{equation}\label{eqexp}
\hat H(\l)=(\l-c\hat H(\l+\d))^{-1}.
\end{equation}
>From this equation, we shall deduce the following 
\begin{prop}\label{Bessel}
Assume that $k(u)=ce^{-\d u}$, for some positive constants
$c$ and $\d$, and let $\hat H(\lambda)$ be the Laplace
transform of the unique solution of (\ref{kr}). Then
\begin{equation}\label{Exact}
\hat H(\lambda)=\frac{J_{\lambda\d^{-1}}(z)}
{J_{\lambda\d^{-1}-1}(z)},\ \ \l>\l_c(\d),
\end{equation}
where $J_\nu(z)$ denotes the Bessel function of order $\nu$ and
$z:=2\sqrt{c}/\d$. Let $j_\nu$ denotes the smallest real positive 
root of $J_\nu$.
Then  $\l_c(\d)$ is given by the equation
$$j_{\l_c(\d)/\d-1}=z.$$
$\l_c(\d)$ is right continuous in $\d$ at zero and satisfies
$$\l_c(\d)=2\sqrt{c}-a  c^{1\over 3} \d^{2\over 3} +O(\delta)$$
with $a\simeq 2,34$. 
\end{prop}
\prf
First notice that $c$ can be chosen equal
to one up to replace $\hat H$ by $\sqrt{c}\hat H(\sqrt{c}.)$
, $\d$ by $\gamma=\d/\sqrt c$ and therefore
$\l_c$
by  $(\l_c/\sqrt{c})$.
Moreover, it is known that if $J_\nu(z)$
is the Bessel function,  and $$h(\nu,z)={J_{\nu}(z)\over J_{\nu- 1}(z)}$$
then
\begin{equation}\label{eqexp2}
h(\nu,z)={{1\over 2}{z\over \nu}\over 1-{z\over 2\nu}h(\nu+1,z)}
\end{equation}
(see \cite{Watson}, chap. 5.6, p. 153). Furthermore, 
$h(\nu,z)$ is uniquely determined by (\ref{eqexp2})
and the boundary condition $\lim_{\nu\ra\infty}h(\nu,z)=0$.
Putting $z={2\over\gamma}$ and $\nu={\l\over \gamma}$,
we find that  since $\hat H$
satisfies (\ref{eqexp}) with the same boundary condition than
$h$, they are related by
$$h(\nu,1)=\hat H(\gamma\nu,1).$$
Therefore, the critical point $\l_c$
corresponds to the largest  $\l$ such that
$$J_{{\l\over\gamma}-1}({2\over\gamma})=0.$$
Since the zeros $j_{\nu,s}$ of the Bessel function
increases with $\nu$ (see 9.5.2 in \cite{AbramoSte})
it follows
that if $j_\nu=j_{\nu,1}$ denotes the smallest zero of
the Bessel function $J_\nu$, the equation for the critical point is
$$j_{{\l_c(\gamma)\over \gamma}-1}={2\over\gamma}.$$
It is known (see 9.5.14 of \cite{AbramoSte}
) that
as $\nu$ is large,
$$j_{\nu}\approx\nu+1,85575 \nu^{1\over 3} +O(\nu^{-{1\over 3}})$$
(for the derivation of this asymptotics, see \cite{Watson},
 chap. XV, 15.83, p. 521, Sturm's method)
so that we deduce (recall that $\nu={\l\over \gamma}$,
$\gamma={\d\over\sqrt{c}}$),
$$\l_c(\d)=2-1,85575\l_c(\d)^{1\over 3}
 \gamma^{2\over 3} +O(\gamma)=2-2.34 \gamma^{2\over 3} +O(\gamma) .$$

\hfill\endproof

\begin{lem}
Let $j_{\nu,n}$, $n\ge 1$, $\nu\in \R^+$,
 be the real positive
zeros of $J_\nu$ arranged in increasing order.
Set $\nu+1=\lambda\d^{-1}$ and $z={2\sqrt{c}\over \d}$. Then
$$\hat H(\lambda)=\frac{2z}{j_\nu^2-j_{\nu_c}^2}+
2z\sum_{n\ge 2}\frac{1}{j_{\nu,n}^2-j_{\nu_c,1}^2},\ \l>\l_c(\d),$$
with
$$\hat H(\lambda)\approx \frac{z J_{\nu_c+1}^2(z)}{2\nu_c \int_0^z
J_{\nu_c}^2(t){\rm d}t/t}\ \ \frac{1}{(\nu-\nu_c)},\ 
\lambda\to\lambda_c(\d),$$
if $j_{\nu}=j_{\nu,1}$.
Hence, 
\begin{equation}\label{limexp}
\lim_{x\ra \infty}\exp(-\lambda_c(\d) x)
H(x)= \frac{z J_{\nu_c+1}^2(z)}{2\nu_c \int_0^z
J_{\nu_c}^2(t){\rm d}t/t}.
\end{equation}
\end{lem}

\prf
The  first identity is a classical result
(see e.g. \cite{Er}, vol.2, p.61). The asymptotic
behavior of $\hat H$ when $\l\to\l_c(\d)$
is obtained by considering the first term
$(2z)/(j_\nu^2-j_{\nu_c}^2)$ and using the
analyticity of the
 smallest positive zero $j_{\nu,1}$ of the
 Bessel function
when the argument is the order $\nu$, using
the asymptotics $(j_\nu^2-j_{\nu_c}^2)
\sim 2j_{\nu_c}
(\partial j_\nu/\partial\nu)_{\nu=\nu_c}(\nu-\nu_c)$,
and standard formulas for the derivative
$\partial j_\nu/\partial\nu$ (see \cite{Watson}).
It remains to consider the problem of asymptotic behavior
of the argument of the Laplace transform.
In our situation,
let 
$$f(s)=\hat H(\l_c(\d)+s)=\int_0^\infty \exp(-su)G(u){\rm d}u,$$
where we set $G(u)=\exp(-\l_c(\d)u)H(u)$. Then
$$f(s)=\frac{{\rm J}_{(\l_c(\d)+s)/\d}(z)}{{\rm J}_{(\l_c(\d)+s)/\d-1}(z)},$$
and the main problem is to find an analytic continuation. Coulomb \cite{Cou}
proved that the roots of the equation in $\nu$,
${\rm J}_\nu(z)=0$ are contained in the real axis when
$z$ is real positive, and the analytic continuation is simply
given by the ratio of Bessel functions where the order $\nu$ is
restricted to $\C\setminus\{(-\infty,\nu_c)\}$. Hence, $f$ can be
continued analytically to
$\C\setminus\{(-\infty,\nu_c)\}$. Further, 
$f$ goes to zero as $|z|\ra\infty$ in
this domain  so that Lemma \ref{Inverse}
applies, yielding (\ref{limexp}). This follows from
classical asymptotics: one uses the expansion
$$J_\nu(z)=\sum_{m\ge 0}\frac{(-1)^m (z/2)^{\nu+2m}}{m!\Gamma(\nu+m+1)},$$
which is analytic as function of $\nu\in\C$ (see e.g. \cite{WatsonB}),
and the Stirling's series for the Gamma function $\Gamma(w)$ when
the complex argument $w\in\C$ is such that
$\vert {\rm arg}(w)\vert\le\pi-\triangle$, for some
positive number $0<\triangle<\pi$
to get
that
$$J_\nu(z)\sim \exp(\nu+\nu\ln(z/2)-(\nu+1/2)\ln(\nu))\sqrt{2\pi},$$
as $\vert \nu\vert\to\infty$ 
(\cite{Watson}, chap. 3).

\hfill\endproof

\subsubsection{Vanishing covariances ; the (almost) general case}

\begin{theo}\label{vancov}
Assume that the kernel $k$ is
 such that there exists $a>1$ and $C<\infty$
so that
$$0\le k(u)\le {C\over (1+x)^{a}}.$$
Then,
$$
\exp(-\l_c(H)t)H(t)\sim \frac{1}{2 A},\ \ x\to +\infty,$$
where
$$A=\frac{{\rm d}}{{\rm d}\l}(\l-\hat{Hk}(\l))\vert_{\l=\l_c(H)}>0.$$
\end{theo}
\prf Hereafter we denote $\l_c$ for $\l_c(H)$.
The proof  goes as follows
\begin{enumerate}
\item We first show that we must have $\l_c=\hat{Hk}(\l_c)$.
This already entails that $$\hat H(\l)
\sim (A(\l-\l_c))^{-1},\quad \l\ra\l_c$$
and hence by Tauberian theorem (see Theorem 4.3 of
chap. V in \cite{Wi}),
$$\int_0^x e^{-\l_c(H) u} H(u) {\rm d}u \sim {1\over A\Gamma(2)}x\quad
x\ra \infty.$$
\item To suppress the integral above,
one has to use in general complex analysis,
Tauberian theorems being then
only valid under additional
monotony properties which are not a priori
satisfied here. To this end, we construct an
analytic continuation $h$ of $\hat H$ in a set of
the form 
\begin{equation}\label{set}
\Gamma_{r,\theta}=S_\theta
\backslash {\rm B}(\l_c,r)
\end{equation}
with
$$S_\theta=\{z\in\C :|\arg(z-\l_c)|\le {\pi\over 2}+\theta\}
\mbox{ and }
{\rm B}(\l_c,r)=\{z\in\C : |z|\le r\}$$
for any $r>0$ and $\theta$ small 
enough.
We show that $h(z)$ goes to
zero
as $|z|$ goes to
infinity and
further
$$\lim_{z\in S_\theta}\sup_{|z-\l_c|\ra 0}
|zh(z)-A^{-1}|=0.$$
We can therefore apply Lemma \ref{Inverse}
and conclude.
\end{enumerate}
To prove the first point,
we proceed by contradiction
assuming that $\l_c>\hat{Hk}(\l_c)$,
and constructing then an analytic continuation
of $\hat H$ in a neighborhood
of $\l_c$, which is a contradiction
with the definition of $\l_c$.
To do that, let us notice that under our
hypothesis, $\hat H$ is bounded on $\l\ge \l_c$,
and therefore $\hat G(\l)=\hat{Hk}(\l)$
is bounded and continuously differentiable since
we assumed $a\ge 1$ (note that $|\hat G'(\l)|\le C\hat H(\l)$). Consequently, $\hat H$ is
also continuously differentiable.
Proceeding by induction, we see that $\hat H$
and $\hat G$ are $\Ca^\infty$ 
at $\l_c$. Let $h_n= n! \hat H^{(n)}(\l_c)$
and $g_n=n! \hat G^{(n)}(\l_c)$.  We now
bound the $h_n$ and $g_n$
by using the idea of majoring
sequences
following Cartan \cite{cartan}, chapter VII.
  Remark
that $V(x,y)= (x-y)^{-1}$ is analytic in a neighborhood
of $(\l_c, g_0)$ such as $U:=\{x: |x-\l_c|\le 3^{-1} |g_0-\l_c|\}
\ts \{y:|y-g_0|\le 3^{-1} |g_0-\l_c|\}$ 
since we assumed $g_0\neq \l_c$.
Therefore, there exists $r>0$ ($r$ can
be taken equal to $3^{-1} |g_0-\l_c|$
according to the choice of the above neighborhood)
 such that for $|x-\l_c|<r$
and $|y-g_0|<r$, 
$V(x,y)=\sum_{p,q} c_{p,q} (x-\l_c)^p(y-g_0)^q$
with,  by Cauchy formula, a finite constant $M$ ($=\sup_{|x-\l_c|<r,
|y-g_0|<r} |V(x,y)|$)
such that 
$$|c_{p,q}|\le {M\over r^{p+q}}.$$
Then, if we consider the formal series $h(\l)=\sum h_n (\l-\l_c)^n$
and $g(\l)=\sum_n g_n(\l-\l_c)^n$, they are 
 formal solutions of
$$\l h(\l)=1+h(\l)g(\l)\qquad\Leftrightarrow \qquad h(\l)=V(\l,g(\l))$$
so that we find that 
for all $n\ge 0$,
\begin{equation}\label{lk}
h_n=P_{n}(g_1,\cdots,g_n; c_{p,q})
\end{equation}
with polynomial functions $P_n$ with non negative
coefficients.  Because the
polynomial functions $P_n$ have
non negative coefficients, we deduce that 
$$|h_n|\le P_{n}(|g_1|,\cdots,|g_n|; |c_{p,q}|).$$
Note that for all $n\ge 1$, since $a\ge 1$,
\begin{eqnarray*}
|g_{n+1}|&=&((n+1)!)^{-1} \int_0^\infty e^{-\l_c u} H(u) k(u) u^{n+1} {\rm d}u\\
&\le& C(n+1)! \int_0^\infty e^{-\l_c u} H(u)  u^{n} 
{\rm d}u =C(n+1)^{-1} |h_n|
\end{eqnarray*}
so that we deduce for $n\ge 1$ 
\begin{equation}\label{jkl2}
n |g_{n}|\le  C P_{n}(|g_0|,\cdots,|g_{n-1}|; {M\over r^{p+q}})
\end{equation}
Now, we can construct a majoring sequence 
by considering the solution of
$$\partial_\l k(\l)= F(\l,k(\l))$$
with $F(x,y)=CM(1-{x-\l_c\over r})^{-1}(1-{y-g_0\over r})^{-1}$ 
and $k(\l_c)=k_0=g_0$.
In fact, by the implicit function theorem
the solution $k$ exists and
is unique in a neighborhood
of $(\l_c)$ such as $U$. 
Writing $k(\l)=g_0+\sum_{n\ge 1} k_n (\l-\l_c)^n$,
we see that $k_i\ge 0$ for all $i$ and 
$$(n+1) k_{n+1}=CP_{n}(k_0,\cdots,k_n; {M\over r^{p+q}})$$
showing with (\ref{jkl2}) by induction 
that for all $n\in\NN$,
\begin{equation}\label{jkl3}
|g_{n}|\le  k_{n}.\end{equation}
Finally, it is not hard
to see that
$$k(\l)=g_0+r\left(1-\sqrt{1+2CM\log(1-{\l-\l_c\over r})}\right)$$
implying that
$$k_{n}\le  \left(r(1-e^{-{1\over 2M}})\right)^{-n}$$
for some finite constant $C$.
This conclude the proof since (\ref{jkl2})
shows that 
\begin{equation}
\label{est}
|g_n|\le \left(r(1-e^{-{1\over 2CM}}) \right)^{-n},
\end{equation} 
so that $g$ is an analytic continuation of $\hat G$
in $\{|\l-\l_c|<r(1-e^{-{1\over 2CM}})\}$ 
and therefore $h(\l)=(\l-g(\l))^{-1}$
is an analytic continuation of $\hat H$ 
in $\{ |\l-\l_c|<2^{-1}r(1-e^{-{1\over 2CM}})\}$. 
This contradicts the definition
of $\l_c$. Thus $\l_c=\hat{G}(\l_c).$

\bigskip
\nn
We now construct an analytic continuation of $\hat H$.
\begin{itemize}
\item Analytic continuation to $S_\theta\cap\{|\Im z|\ge R\}$
for some sufficiently large $R$ :
>From $\l_c=\hat{G}(\l_c)$ and the fact that $\hat G$
is continuously differentiable
at $\l_c$, we see that 
$$\hat H(\l)\sim (A(\l-\l_c))^{-1}\quad \l\ra\l_c$$
implying by Tauberian theorem (see \cite{Wi})
that 
$$\int_0^x e^{-\l_c t} H(t) {\rm d}t \sim {1\over \Gamma(2) A} x,
\quad  x\ra\infty.$$
Consequently, we see by integration
by parts that if $a>1$,
$$B:=|\hat{Hk}(\l_c)|<\infty$$
and therefore, by (\ref{f5}),
that the analytic continuation  of
$\hat H$ to $\{\Re(z)>\l_c\}$ satisfies uniformly
on this set
$$\hat H(z)=z^{-1} +o(|z|^{-1}),\quad |z|\ra \infty$$
In particular, for $R$ large enough,
$z-\hat G(z)$ has no pole 
in $\Gamma_R=\{\Re(z)>\l_c\}\cap {\rm B}(\l_c,R)^c$.
We can therefore proceed
as above by majoring sequences
to see that we can extend analytically
$\hat H$  around each point $z_0$
of the type $z_0=\l_c+\varepsilon +i y$
with $y>R$ and $\varepsilon>0$
and this continuation
is analytic in $|z-z_0|<C |y-R|$
for some universal constant $C$ (indeed note that
here the radius  $r$ of
convergence of $V$ is of the order of the distance $|y-R|$). 
Further, for $y\ge 2R$
it is not hard to
see that on $|z-z_0|<C |y-R|$, the continuation of $\hat H$ 
and therefore $\hat G$,
is bounded by $B$. Performing such analytic
continuation for every $z_0=\l_c+\varepsilon +i y$
with $|y|\ge R$, we obtain an analytic continuation
$h,g$ of $(\hat H,\hat G)$ on $\Gamma_R=S_\theta\cap 
 \{|\Im (z)|\ge R\}$ 
with $\tan(\theta)\le 2^{-1}C$ such that $g$ is uniformly
bounded by $B$.
Moreover, note that since by construction, $\hat G$ remains 
uniformly bounded and the continuation satisfies $\hat H=(z-\hat
G)^{-1}$,
$$\hat H(z)=O(|z|^{-1}),\quad
|z|\ra\infty$$

\item  Analytic continuation to $\{|\arg(z-\l_c)|\le{\pi\over
2}+\theta\}\cap {\rm B}(\l_c,r)^c \cap {\rm B}(\l_c,R)$:
Again, the main issue 
is to control the zeros of $z-\hat{Hk}(z)$.
Let us study these zeroes on $\Re(z)=\l_c$.
Observe that for such a $z$,
$$z-\hat{Hk}(z)=0=\l_c-\hat{Hk}(\l_c)$$
Taking the real part of both sides of
this equality, we find that
$$\int_0^\infty H(u) k(u) e^{-\l_c u} (\cos(\Im(z)u)-1) {\rm d}u=0$$
which implies that $\Im(z)=0$ since $Hk\ge 0$.
Hence, $(\l_c,0)$ is the only zero on $\Re(z)=\l_c$.
We can again apply majoring sequences 
to continue $\hat H$ in the neighborhood
of any $(\l_c, y)$ with $y\neq 0$
in some domain ${\rm B}((\l_c, y),r_y)$
for some $r_y>0$. We thus obtain a continuation
on $\cup_{r\le |y|\le R} {\rm B}((\l_c, y),r_y)
$ which can be reduced to
a finite union $\cup_{1\le i\le L} {\rm B}((\l_c, y_i),r_{y_i})
$ since we are in the compact $B(\l_c,R)$.
Thus, since $\e =\min_{1\le i\le L} r_{y_i}>0$,
we obtain an analytic continuation of
$\hat H,\hat G$ on $$\{\Re(z)>\l_c-\e\}\cap 
{\rm B}(\l_c,r)^c \cap {\rm B}(\l_c,R)\subset S_\theta
\cap 
{\rm B}(\l_c,r)^c \cap {\rm B}(\l_c,R) $$
where the latter inclusion holds
for $\theta\le \theta_0\sim  {\e\over R}$.
This finishes the construction
of the continuation 
of  $(\hat H,\hat G)$.

\end{itemize}
Finally, noting that for $z\in S_\theta$
approaching $\l_c$, the differentiability 
 of $\hat G$ shows that
$z-\hat G(z)\sim A(z-\l_c)$
with $A=1-G'(\l_c)\ge 1>0$ since $G$ is decreasing, 
which implies 
$$\hat H(z)\sim (A(z-\l_c))^{-1}.$$
We can thus conclude the proof of the lemma thanks 
to Lemma \ref{Inverse}.

\hfill\endproof


\subsection{Covariances with non zero
limit}
In this section we consider the case where
$$\lim_{t\ra\infty} k(t)=C> 0.$$
We first tackle the case 
where the covariance decays towards this limit 
exponentially fast (which
is somewhat simpler) and then when the speed is only algebraic.

\subsubsection{Exponentially decaying covariances\label{ExpoDecay}}

Let us now assume that we have $k(u)=c_2+c_1e^{-\d u}$ with $c_1,c_2>0$,
which should correspond to the case where we consider
the
 p-SSK model in the range $t/s$ of order one
 (see the introduction or \cite{CuKu}). Then, we obtain
$$\l \hat H(\l)=1+\hat H(\l)(c_2 \hat H(\l) +c_1\hat H(\l+\d)),\ \ \l>\l_c(H).$$
We can solve this equation to find that for $\l>\l_c(H)$
\begin{equation}\label{f7}
\hat H(\l)=(2c_2)^{-1}[\l -c_1\hat H(\l+\d)-\sqrt{(\l-c_1\hat
H(\l+\d))^2-4c_2}],
\end{equation}
implying  $\l_c-c_1\hat H(\l_c+\d)\ge 2\sqrt{c_2}$.
Set $\l_c:=\l_c(H)$.

We claim
that $\l_c$ is the unique positive  number
such that $$\l_c-c_1\hat H(\l_c+\d)=2\sqrt{c_2}.$$
Suppose that $\l_c$ is such that
$\l_c-c_1\hat H(\l_c+\d)>2\sqrt{c_2}$. Using the 
analyticity of $\hat H$ on its domain of convergence,
this remains true for $\l$ with
$\vert \l-\l_c\vert<\varepsilon$,
for some small enough positive constant
$\varepsilon<\d$. Let
$D:=\{(x,y)\in\R^2;\ x-c_1 y>2\sqrt{c_2}\}\subset\R^2$, and
consider the mapping $\Psi:\ D\longrightarrow \R$ given by
$\Psi(x,y):=(2c_2)^{-1}[x-c_1 y
-\sqrt{(x-c_1y)^2-4c_2}]$.
Then, the function $h:\ (\l_c-\varepsilon,+\infty)
\longrightarrow \R$
given by $h(\l):=\Psi(\l,\hat H(\l+\d))$ is analytic
with $h(\l)=\hat H(\l)$, $\forall \l>\l_c$, and
provides thus an analytic continuation of $\hat H$
on a domain containing
its domain of convergence, a contradiction
with the fact that the abscissa of convergence
$\l_c$ is a singularity of $\hat H$ when
$H$ is non-negative (see Theorem 5b in \cite{Wi}).
The computation of the Laplace transform seems
difficult to obtain in closed form, and the
abscissa of convergence $\l_c$ remains
unknown.
\begin{theo}\label{secondordermixed}
Assume that $k(u)=c_2+c_1\exp(-\d u)$, for
positive constants $c_1$, $c_2$ and $\d$.
Then
$$\hat H(\l)=\frac{1}{\sqrt{c_2}}-\frac{1}{\sqrt{c_2}}(\l-\l_c)^{1/2}+o((\l-\l_c)^{1/2}),$$
and
$$\exp(-\l_c x)H(x)\sim A x^{-3/2},\ \ x\to+\infty,$$
for some positive constant $A$.
\end{theo}
\prf
 First note that
$$\lim_{\l\to\l_c}\hat H(\l)=\lim_{\l\to\l_c}h(\l)
=(2c_2)^{-1}(\l_c-c_1\hat H(\l_c+\d))=\frac{1}{\sqrt{c_2}}.$$
Next,
$\hat H(\l)-(2c_2)^{-1}(\l_c-c_1\hat H(\l_c+\d))$ is given by
$$-\frac{1}{2c_2}\sqrt{(\l-c_1\hat H(\l+\d))^2-4c_2},$$
with
$$\l-c_1\hat H(\l+\d)-2\sqrt{c_2}=(1-c_1\hat H(\l_c+\d)')(\l-\l_c)+o
(\l-\l_c),$$
and
$$0<
1-c_1\hat H(\l_c+\d)'=
1+c_1\int_0^{+\infty}t e^{-(\l_c+\d)t}H(t){\rm d}t<+\infty.$$
Thus,
$$\hat H(\l)\approx \frac{1}{\sqrt{c_2}}-\frac{1}{c_2^{3\over 4}}
(1-c_1\hat H(\l_c+\d)')^{\frac{1}{2}} (\l-\l_c)^{1/2}+o((\l-\l_c)^{1/2}),$$
Moreover,
when $\l>\l_c$, 
$$\hat H'(\l)=
\frac{(1-c_1\hat H'(\l+\d))}{2c_2}(1-
\frac{(\l-c_1\hat H(\l+\d))}{\sqrt{(\l-c_1\hat H(\l+\d)^2-4c_2}}),$$
with
\begin{equation}\label{eqh}
\sim -\frac{(1-c_1\hat H'(\l_c+\d))^{1/2}}{(2c_2)^{3/4}}
(\l-\l_c)^{-1/2},\ \ \l\to\l_c.
\end{equation}
To prove the theorem, we need as 
for the proof of Theorem \ref{vancov},
to continue $\hat H$ (and therefore $\hat H'$) analytically on sets
of the form (\ref{set}). We in fact continue
it on $[S_\theta\cap B(\l_c,R)\backslash B(\l_c,r)]
\cup \{|\Im(z)|\ge R\}$ for some sufficiently large $R$ and small 
$\theta$.
Let 
$$\triangle_z=(z-\d-c_1\hat H(z))^2-4c_2$$
 with
$\triangle_{\l_c+\d}=0$ (see above). First notice that
$$\vert \hat H(z)\vert\le\int_0^\infty \exp(-\Re(z)x)H(x){\rm d}x
:=A,$$
so that $\hat H$ is uniformly bounded on $\{\Re(z)\ge\l_c\}$.
Therefore, for $R$ large and any $z$ such
that $\Re(z)>\l_c-\d$
and $z\in {\rm B}(\l_c,R)^c=\{z:|z-\l_c|\ge R\}$, 
$\D_{z+\d}\approx z^2+O(1)\ne 0$. 
Hence, if we set, for $\{\Re(z)\ge \l_c-{1\over 2}\d\}\cap
\{|\Im(z)|\ge R\}$
$$h_1(z)=(2c_2)^{-1}[z-c_1 H(z+\d)-\sqrt{\D_{z+\d}}$$
is analytic and thus provides an analytic continuation
of $\hat H$. Further,
note that if $R$ is large enough,
$$\sup_{ z\in B(\l_c,R)^c\cap \{\Re(z)\ge \l_c-{1\over 2}\d\} }
 |h_1(z) |\le \sup_{ z \in B(\l_c,R)^c\cap\{\Re(z)\ge\l_c\}}
|\hat H(z) |:=A.$$
Indeed, 
\begin{eqnarray*}
|h_1(z) |&=&|z-c_1\hat H(z+\d)|^{-1}|1+(1-4c_2(z-
c_1\hat H(z+\d))^{-1})^{1\over 2}|^{-1}\\
&\le &(|z|-c_1A)^{-1}
|1+(1-4c_2(|z|-c_1A)^{-1})^{1\over 2}|^{-1}\le A,
\end{eqnarray*}
where the last inequality holds for $R$ large enough.
>From this formula, we may proceed by induction to
construct an analytic continuation
of $\hat H$ on $|\Im(z)|\ge R$
by arguing by induction that $\D_{z+\d}$ does not vanish.
This continuation remains uniformly bounded 
by $A$.

We next show that $\D_{z+\d}\neq 0$ for $z\in {\rm B}(\l_c,R)\cap
  S_\theta$ for $\theta$
small enough. Indeed, 
the analyticity of $\hat H$ on its domain
implies then that the compact intersection $D_{u,R}$
of ${\rm B}(\l_c,R)$ with $\{z\in\C;\ \Re(z)\ge\l_c+u\}$,
$u>0$, contains only a finite number of roots
of the equation $\triangle_z=0$, and therefore  
$\D_{z+\d}$ has only finitely 
many roots in ${\rm B}(\l_c,R)\cap
  S_\theta$. We can thus choose $\theta_0>0$
such that $\D_{z+\d}\neq 0$ for $z\in {\rm B}(\l_c,R)\cap
  S_\theta$ when  $\theta\le\theta_0$. As a consequence, if we let 
for  $\varepsilon >0$ small enough,  the domain
$\Gamma_\varepsilon^1$ be  given by
$$\Gamma_\varepsilon^1=\{z\in\C\setminus
 \{\l_c\};\ \l_c-\d+\varepsilon<\Re(z)<\l_c+
\varepsilon\}\cap S_{\theta_0},$$
and define the function $\Psi(z)$ on this domain as
$$\Psi(z)=(2c_2)^{-1}[z-c_1\hat H(z+\d)-\sqrt{\triangle_{z+\d}}].$$
Then $\Psi$ is analytic and, from (\ref{f7}),
coincides with $\hat H$ on the band $\{z\in\C;\ \l_c<\Re(z)<\l_c+\varepsilon\}$,
and provides thus an analytic continuation of $\hat H$
on $\Xi_\varepsilon:= \Gamma_\varepsilon^1\cup \{\Re(z)\ge\l_c\}$.

At the end of the day, we have constructed
an analytic continuation of $\hat H$ on $S_{\theta_0}$.
Further, because it remains uniformly bounded,
we also see that for large $|z|$,
$$|\hat H(z)|\approx O({1\over |z|}).$$
Consequently, $\hat H'$ can also be extended analytically to
$S_{\theta_0}$ and its continuation
 remains uniformly bounded too. 
As a consequence,
$$|\hat H(z)|\approx O({1\over |z|})$$
and $|\hat H'(z)|\approx O({1\over |z|^2})$
by (\ref{eqh})
We can therefore apply 
 Lemma \ref{Inverse} and conclude.

\hfill \endproof

\subsubsection{Algebraically decaying covariances\label{AlgeDecay}}

We consider here the case where
the stationary covariance takes the form
$k(u)=c_2+c_1 k_1(u)$, for some
positive constants $c_1$, $c_2$ and $|k_1(u)|\le (1+u)^{-1})$
for some $a\ge 1$.
Set for convenience $\l_c=\l_c(H)$,
with $\l_c=\l_c(Hk)$. The basic relation becomes, for $\l>\l_c$,
\begin{equation}\label{f77}
\l\hat H(\l)=1+c_2\hat H(\l)^2
+c_1\hat G(\l)\hat H(\l),
\end{equation}
where we set
$$\hat G(\l)=\int_0^\infty 
\exp(-\l u)H(u)k_1(u){\rm d}u,$$
which converges for $\l>\l_c$, the Laplace transform of the 
function $G(u)=H(u)k_1(u)$.

We shall prove that

\begin{theo}\label{algdec}

\begin{enumerate}
\item $\l_c$ is solution of the equation
 $$\l_c-c_1\hat G(\l_c)=2\sqrt{c_2}.$$
\item 
$$\hat H(\l)=\frac{1}{\sqrt{c_2}}-\frac{1}{\sqrt{c_2}}(\l-\l_c)^{1/2}+o((\l-\l_c)^{1/2}),$$
and
$$e^{-\l_c t} H(t)\sim A t^{-3/2},\ \ t\to+\infty,$$
for some positive constant $A$.
\end{enumerate}
\end{theo}
\prf 
Note that $\hat H$ is uniformly
bounded on $ \l\ge\l_c$
by $\sqrt{c_2}^{-1}$
so that the integral 
defining $\hat G$ is 
absolutely
convergent and  $\lim_{\l\to\l_c}\hat G(\l)$ exists. 
(\ref{f77}) also gives
the equation in $\hat H$, $c_2\hat H^2+(c_1\hat G-\l)\hat H+1=0$,
$\l>\l_c$, showing that the discriminant
$\triangle_\l=(\l-c_1\hat G)^2-4c_2$ is non-negative.
Thus
\begin{equation}\label{in1}
\l_c-c_1\hat G(\l_c)\ge 2\sqrt{c_2}
\end{equation}
and  for $\l>\l_c$,
\begin{equation}\label{in2}
\hat H(\l)=\frac{1}{2c_2}(\l-c_1\hat G(\l)-\sqrt{\triangle_l})
,\ \ \l>\l_c,\end{equation}
where the branch was chosen to satisfy the condition 
 $\lim_{\l\to +\infty}\hat H(\l)=0$
and $\triangle_l=(\l-c_1\hat G(\l))^2-4c_2$.

We can proceed exactly as in the proof of theorem \ref{vancov}
;
we prove that (\ref{in1})
is an equality by contradiction
using majoring sequences. The analytic continuation is
also obtained similarly.

\hfill\endproof

\section{More general limiting behaviour}\label{general}

Let us consider the case where
 $k$ has
a 
sufficiently
flat part around the diagonal and a stationary part.
Consider
a   non negative function $h$  such that,
there exists a positive constant $C$
and for $T,M>0$ a function $\d(T,M)$ such that
$\d(T,M)\ra 0$ for all $M$ when $T$ goes to
infinity
so that

\begin{equation}\label{hyp}
\sup_{t\ge T}\sup_{|s-t|\le M}| h(s,t)-C|\le
\d(M,T),\quad \sup_{s\ge t\ge T}|h(s,t)|\le C.
\end{equation}
Note that the second condition is a consequence of the
first
when $h$ is a covariance, which we shall not need
to assume.
For instance, it is clear that
such an assumption is verified by
the ratio
$$ h(s,t)
=C \left({t\over s}\right)^a \quad\mbox{ for }t<s,$$
with some $a\ge 0$ or any linear combination of
such functions.

Then, we claim that

\begin{theo}\label{flat}
Let $k$ be a covariance kernel
such that
$$k(s,t)= k_1(s-t)+h(s,t)$$
with $h$ satisfying (\ref{hyp}) and 
$k_1$ is a non negative function.
Then,
regardless of
the way $t$ goes to
infinity 
$$\lim_{t\ra\infty} \lim_{s-t\ra\infty}{1\over s-t}\log
H_k(s,t)=\lim_{t\ra\infty}\lim_{s-t\ra\infty}{1\over s-t}\log
H_{C+k_1} (s,t)=\l_c(H_{C+k_1}).$$
\end{theo}
Hence, this theorem shows that the first order asymptotics of $H$ are
only governed by its stationary
part. As a direct consequence,
\begin{cor}
Let $h$ satisfying (\ref{hyp}).
Regardless of
the way $t\ge T$ goes to
infinity, if $h\ge 0$,
$$\lim_{s-t\ra\infty}{1\over s-t}\log
H_{h}(s,t)=2\sqrt{C}$$
and
$$\lim_{s-t\ra\infty}{1\over s-t}\log
H_{h+C'e^{-\d|s-t|}}(s,t)=\l_c(H_{C+C'e^{-\d .}}).$$
\end{cor}

\nn
{\bf Proof of theorem \ref{flat} :}
By property \ref{increas}
and the second hypothesis 
in (\ref{hyp})
$$H_k(s,t)\le  H_{C+k_1}(s,t)$$
resulting with the announced upper bound.
For the lower bound,
note that since $k$ is non negative,
for any $n\in\N$, any $K\in\N$, 
\begin{eqnarray}
H(s,t)&\ge&
\sum_{\s\in\NC_n} \int_{t\le t_1\cdots
t_{2n}\le s}\prod_{i=1}^n k(t_i,t_{\s(i)})
\prod_{j=1}^{2n} {\rm d}t_j\nonumber\\
&\ge& \sum_{\s\in\NC_n^K}
\int_{t_1,\cdots,t_{2n}\in\D_n^K}
\prod_{i=1}^n k(t_i,t_{\s(i)})
\prod_{j=1}^{2n} {\rm d}t_j.\label{turlut}\\
\nonumber
\end{eqnarray}
Here,
$\NC_n^K$ are the elements of $\NC_n$
where partitions occur only inside
the boxes $[2Kp, 2K(p+1)]$ for $p\in \{0,\cdots, [{n\over K}
]-1\}$ or $[2K[{n\over K}
], 2n]$. In other words, crossing between these
boxes are prohibited and
$\NC_n^K$ is given by the
set of non-crossing involutions $\sl$
of $\NC_n$ such that
$\s|_{[2Kp, 2K(p+1)]}\in \NC_K$,
$\forall  0\le p\le [{n\over K}]-1$,
and
$ \s|_{[2K[{n\over K}
], 2n]}\in \NC_{n-K[{n\over K}
]}$.
Moreover,
$$\D_n^K=\big\lbc t+(s-t){p K\over n}\le t_{2Kp +1}\le\cdots
\le t_{2K(p +1)}\le t+(s-t){(p+1) K\over n}, 
 \qquad$$
$$\quad\qquad 0\le p\le [{n\over K}]-1,\ t+(s-t){K\over n}
[{n\over K}
]\le  t_{2(n-K[{n\over K}
])+1} \cdots\le t_{2n}\le s\big\rbc.$$
Observe that by construction,
when $\s\in \NC_n^K$, for all $i$,
$t_i$ and $t_{\s(i)}$ belong to the same box
of the partition $\D_n^K$.
Hence by our hypothesis, for $t\ge T$,
for all $\s\in \NC_n^K$, all ${\bf t}\in \D_n^K$,
all $i\in\{1,\cdots,n\}$, 
$$ k(t_i,t_{\s_i})\ge k_1(t_{\s_i}-t_i)+
 \inf_{|t'-s'|\le {K\over n}
(s-t)} h(t',s')\ge k_1(t_{\s_i}-t_i)+C-\d$$
provided $\d( {K\over n}
(s-t),t)\le\d$. 
Therefore, we deduce from
(\ref{turlut}) that
$H(s,t)$ is larger than
$$
\prod_{p=1}^{[{n\over K}]} 
 \sum_{\s\in\NC_K}  \int_{t+(s-t){p K\over n}\le t_{2Kp +1}\le\cdots
\le t_{2K(p +1)}\le t+(s-t){(p+1)K\over n}}\qquad$$
$$\quad\qquad\prod_{i=1}^{K} (k_1(t_{\s_i}-t_i)+C-\d)
\prod_{i=1}^{2K}{\rm d}t_i$$
$$
\ts (C-\d)^{n-K[{n\over K}]}
C_{n-K[{n\over K}]} {[(s-t)(1-{K\over n}[{n\over K}])]^{2(n-K[{n\over
K}])}\over 2(n-K[{n\over K}])!},
$$
where in the last line we bounded below the term corresponding to the
indices betwen
$2[{n\over K}]K$ and $n$ with the convention $0^0=1$.
It is not hard to see that 
we can neglect this correction term 
(indeed, we shall take later 
$n$ of order $s-t$ 
and $K$ large, but independent of $s-t$).
As a consequence of the above lower bound,
we have that if we define 
$ B^K_{C-\d +k_1}( {(s-t)K\over n})$, $\d>0$, by 
$$
\sum_{\s\in\NC_K}  \int_{t+(s-t){p K\over n}\le t_{2Kp +1}\le\cdots
\le t_{2K(p +1)}\le t+(s-t){(p+1)K\over n}}$$
$$\prod_{i=1}^{K} (k_1(t_{\s_i}-t_i)+C-\d)
\prod_{i=1}^{2K}{\rm d}t_i,$$
then,
for any $K,n, s-t,t$ such that $\d(t, {(s-t)\over n}K)\le \d$
\begin{eqnarray*}
H(s,t)&\ge& [ B^K_{C-\d +k_1}( {(s-t)K\over n})]^{n\over K}\\
&=&[ B^K_{C-\d +k_1}(u)]^{s-t\over u}\\
\end{eqnarray*}
where we have set
$u={s-t\over n}K$.
Now, using Jensen's inequality when
${n\over K}={s-t\over u}>1$, we deduce for any $C'>0$,
$n>K$, $t,s$ so that $\d( {(s-t)\over n}K,t)\le \d$, 
\begin{eqnarray}
H(s,t)&\ge& {1\over 2C'eu} \sum_{K\le 2C'eu}
[ B^K_{C-\d +k_1}(u)]^{s-t\over u}\nonumber\\
&\ge&[{1\over 2C'eu} \sum_{K\le 2C'eu}
 B^K_{C-\d +k_1}(u)]^{s-t\over u}\label{blurp}\\
\nonumber
\end{eqnarray}
Recall that if  $C'=C+||k||_\infty$
we already observed that with $A=\sum n^{-2}$
\begin{equation}\label{plu}
H_{C-\d +k_1}(u)\le A\max_{n\le 2 C' e u} n^2 B^n_{C-\d +k_1}(u)
\le A (2 C' e u)^2 \sum_{n\le 2 C' e u}B^n_{C-\d +k_1}(u).
\end{equation}
Thus, we deduce from (\ref{blurp}) that 

\begin{eqnarray}
H(s,t)
&\ge& [{1\over A(2C'eu)^3 }H_{C-\d +k_1}(u)]^{s-t\over u}\label{klj}\\
\nonumber
\end{eqnarray}
where we have used (\ref{plu})
in the last line.
Now, for any $\d>0$ by   section \ref{stationary},
there exists $\l_c(C-\d +k_1)>0$ such that
$$\lim_{u\ra\infty}{1\over u}\log
H_{C-\d +k_1}(u)=\l_c(C-\d +k_1)$$
so that we arrive at, for any $\e>0$, for $u\ge u(\e)$ 
large enough, ${s-t\over u}>1$, $\d(u,t)<\d$

\begin{eqnarray}
H(s,t)&\ge& [{1\over A(2Ceu)^3 }]^{s-t\over u}
e^{(\l_c(C-\d +k_1)-\e)(s-t)}
\label{klj2}\\
\nonumber
\end{eqnarray}
which shows, by taking first $s-t$ going to
infinity while $u\ge u(\e)$ is fixed,
$$
\liminf_{s-t\ra\infty}
{1\over s-t} \log H(s,t)\ge {1\over u}\log(A(2Ce u)^3)+\l_c(C-\d+k_1)-\e$$
and then letting
 $t$ going to infinity (and hence $u$ going to zero),
and finally $u$ going to infinity, 
$$
\liminf_{s-t\ra\infty}
{1\over s-t} \log H(s,t)\ge\lim_{\d\downarrow 0}
\l_c(C-\d +k_1).$$
Property (\ref{weak})
completes the proof since $k_1\ge 0$ implies $C-\d+k_1\ge (1-\d
C^{-1}) (C+k_1)$.

\hfill \endproof

\nn
{\bf Acknowledgments }
We are very grateful to P. Gerard 
and R. Speicher for
cheerful and motivating discussions.





\end{document}